\documentclass[a4paper,USenglish,cleveref, autoref, thm-restate]{lipics-v2021}
\nolinenumbers
\pdfoutput=1 
\hideLIPIcs  
\graphicspath{{./img/}}
\usepackage[utf8]{inputenc}
\usepackage[T1]{fontenc}
\usepackage{todonotes}
\let\svtodo\todo
\renewcommand\todo[1]{\svtodo[inline]{#1}}

\title{Formalizing a Diophantine Representation of the Set of Prime Numbers
(short paper)}
\titlerunning{Formalizing a Diophantine Representation of the Set of Prime Numbers}

\author{Karol P\k{a}k}{University of Bia\l{}ystok, Poland}{pakkarol@uwb.edu.pl}{https://orcid.org/0000-0002-7099-1669}{} 

\author{Cezary Kaliszyk}{University of Innsbruck, Austria}{cezary.kaliszyk@uibk.ac.at}{https://orcid.org/0000-0002-8273-6059}{}

\authorrunning{K. P\k{a}k and C. Kaliszyk}
\Copyright{Karol Pąk and Cezary Kaliszyk}



\supplement{Formalization can be found at:}
\supplementdetails{Formalization}{http://cl-informatik.uibk.ac.at/cek/itp2022/}
\funding{ERC starting grant no.~714034 \emph{SMART} and Cost action CA20111 \emph{EuroProofNet}}
\acknowledgements{We would like to thank Yuri Matiyasevich for his comments on the
previous version of this paper.}

\EventEditors{June Andronick and Leonardo de Moura}
\EventNoEds{2}
\EventLongTitle{13th International Conference on Interactive Theorem Proving (ITP 2022)}
\EventShortTitle{ITP 2022}
\EventAcronym{ITP}
\EventYear{2022}
\EventDate{August 7--10, 2022}
\EventLocation{Haifa, Israel}
\EventLogo{}
\SeriesVolume{237}
\ArticleNo{27}
\lstdefinelanguage{Mizar}%
{columns=fullflexible,
keywords={scheme,schemes,environ,provided,where,ranks,
theorem,definition,radix,reserve,properties,struct,inhabited,expandable,attribute,
adjective,registration,coherence,defpred,cluster,from,sch,%
given,such,that,%
reflexivity, irreflexivity, symmetry, asymmetry, connectedness,and,attr,%
    antonym,existence,uniqueness,commutativity,idempotence,synonym,notation,%
    mode, means,func,pred, pred,equals,it,of,is,axiomatization,sethood,reconsider,redefine,%
    if,otherwise,proof,for,ex,being,holds,def,let,consider,take,not,contradiction,st,the,%
    be,thus,implies,assume,then,not,by,or,hence,thesis,end,iff,;,:,",",\#,as,qua},%
   sensitive=true,%
   basicstyle={{\linespread{1.}\usefont{T1}{lmss}{m}{n}}},%
 keywordstyle={\usefont{T1}{lmss}{sbc}{n}\selectfont},%
 keywordstyle=[2]{\it},%
  mathescape = true,%
   morecomment=[l][\texttt]{::},%
   literate={'}{{$\strut\mkern6mu\strut$}}1%
   {&}{{\usefont{T1}{lmss}{sbc}{n}\selectfont\texttt{\&}}}1%
   {-}{{\usefont{T1}{lmss}{m}{n}\selectfont\texttt{-}}}1%
   {->}{{$\rightarrow$}}1%
   {BrLeft}{{$\{$}}1%
   {\{\}}{{\{\!\}}}1
  }

\makeatletter

\def\mizV[#1]{\lstset{language=Mizar,keywords=[2]{#1}}\lstinline}
\makeatother
\lstnewenvironment{Mizar}[1]{\lstset{language=Mizar,keywords=[2]{#1}}}{}

\begin{document}

\maketitle

\begin{abstract}
The DPRM (Davis-Putnam-Robinson-Matiyasevich) theorem is the main
step in the negative resolution of Hilbert's 10th problem. Almost
three decades of work on the problem have resulted in several
equally surprising results. These include the existence of
diophantine equations with a reduced number of variables,
as well as the explicit construction of polynomials that represent
specific sets, in particular the set of primes.
In this work, we formalize these constructions in the Mizar system.
We focus on the set of prime numbers and its explicit representation
using 10 variables. It is the smallest representation known today.
For this, we show that the exponential function is diophantine, together
with the same properties for the binomial coefficient and factorial.
This formalization is the next step in the research on formal
approaches to diophantine sets following the DPRM theorem.

\end{abstract}

\section{Introduction}\label{s:intro}

Hilbert's 10th problem (\texttt{H10}) asks whether there exists an
algorithm\footnote{Today interpreted as an adequate RAM program or
  equivalently a Turing machine searching for solutions.}  that can
determine if a diophantine equation has a solution over the integers.
A major step towards the negative resolution of the problem was
achieved by the \emph{Davis conjecture}, stating that the notions of
diophantine sets and recursively enumerable sets coincide. This is the
case since recursively enumerable sets without algorithms for
recognizing their elements have already been known.
Indeed, by the Davis Normal Form Theorem \cite{DBLP:journals/jsyml/Davis53},
for every recursively enumerable set $R\subseteq \mathbb{N}^m$
there exist a number $n$ together with a polynomial $P$ over
$m+n+2$ variables ($n$ of the variables are parameters and $m+2$ are unknowns) with integer coefficients, such that
\begin{equation}\label{eqDNF}
\forall_{a_1,\ldots,a_n} R(a_1,\ldots,a_n) \iff \exists_x \forall_{y\leq x} \exists_{x_1\leq x,\ldots,x_m\leq x}
P(a_1,\ldots,a_n,x,y,x_1,\ldots,x_m)=0
\end{equation}
Eliminating the single universal quantifies from \autoref{eqDNF}, it becomes
the condition defining a diophantine relation. However, the research undertaken
by
Robinson,
Davis,
and
Putnam
to eliminate this quantifier required almost 30 years.
It is therefore not surprising that they created several theorems that give
a negative solution to the problem but under certain assumptions, which enable
such quantifier elimination.
One of such assumptions is that the exponential function can be defined in a diophantine way.
This has been eliminated by
Matiyasevich, who definitively completed the proof of the DPRM-theorem \cite{DPR1961,YuriOrg}, all
every recursively enumerable set of natural numbers is diophantine.

Our work formalizes multiple consequences and results that originated from the proof of the DPRM,
all concerning diophantine equations. We formalize the fact that the exponential function is diophantine
\cite{hilb10_1}, together with the same property for the binomial coefficient and factorial.
These allow the proof of DPRM-theorem \cite{hilb10_5} in a
post-Matiyasevich approach proposed by Smorynski \cite{LNT:Smorynski},
where the concept of recursively enumerable is defined using \autoref{eqDNF}.

We also formalize Robinson's \cite{JR69} conditions for the DPRM-theorem, namely that the set of
primes is representable by a diophantine polynomial.
This polynomial is sufficient to prove the DPRM-theorem in the Robinson approach.
On the other hand, the existence of this polynomial is guaranteed,
but not its explicit statement. In fact, stating it explicitly has been a long-standing challenge and
the problem of finding the polynomial defining the prime numbers with a minimum
number of variables is an open problem in number theory.

In 1971, Yuri Matiyasevich proposed the construction of a diophantine polynomial of degree 37
with 24 variables {(one of the variables is a parameter and 23 are unknowns)} that defined the set of prime numbers. This result has been improved together
with Robinson \cite{MR75} to a polynomial with 14 variables including 13 unknowns. To do this, they showed that every
diophantine polynomial, can be reduced to 13 unknowns. Wada et al. \cite{AMM76} have later reduced
the polynomial to 12 variables, however, the rank of the polynomial is 13,697.

The main result of the current work is the formalization of the polynomial over 10 variables that
defines the prime numbers together with a large number of formalized results in numeral analysis (263 proved top-level MML theorems totaling 922KB) necessary for this result.
This polynomial, proposed by Matiyasevich \cite{eng_30}, is today the smallest known\footnote{Private email exchange with Yuri Matiyasevich, January 2022.}
polynomial for the open problem. This work improves on the work by the first author, presented at the FMM 2021 workshop \cite{hilb10_6},
where the 26-variable polynomial was defined in Mizar without any further properties.

\section{DPRM Formalizations and Their Relation to Number Theory}

Larchey-Wendling and Forster  \cite{H10coq} formalized the DPRM theorem in Coq and
Bayer et al. \cite{DBLP:conf/itp/BayerDPSS19} in Isabelle/HOL.
Both formalizations develop register machines, Minsky machines, advanced properties of the Pell equation
and prove that exponentiation is diophantine.
The first work proved the bounded universal quantification theorem in order to reach the final DPRM theorem,
and then additionally showed the undecidability of \texttt{H10} and discusses other
undecidable problems as future work.
The second work proves the DPRM theorem
following the approach proposed by Matiyasevich in \cite{Mat2000}
instead of proving
the bounded universal quantification theorem. There, the discussed future work is to
extend it to register machines to prove the undecidability
of the Halting problem.
Carneiro's formalization in Lean \cite{carneiro2018matiyasevic} uses Pell equations to prove the key lemma of
Matiyasevich, stating that exponentiation is diophantine.

Our work focuses on the applications of \texttt{H10} listed by Sun \cite{sun2021results} and instead of
register machines we focus on the theory of polynomials.
By \texttt{H10}, $\exists_{x_1,\ldots,x_\nu\in\mathbb{Z}}P(a,x_1,\ldots,x_\nu)=0$ is
undecidable for some diophantine polynomial $P$, $\nu\in\mathbb{N}$. In 1970 Matiyasevich justified $\nu<200$.
Further ingenious number-theoretic ideas allowed him to prove together with Robinson \cite{MR75},
$\nu\geq 13$ by developing general diophantine polynomial reduction methods.
Observe, that when $\nu\geq 13$, the unknowns range over positive (or non-negative) integers.
In 1975, Matiyasevich further announced that $\nu\geq 9$ and Jones gave a complete proof \cite{Jones82},
but in both cases (Matiyasevich's announcement and Jones's proof) the range of all unknowns is again limited to the positive (or non-negative) integers.
Sun \cite{sun2021results} improved this result by modifying the range of the unknowns.
He showed the case $\nu \leq 9$ with only one limitation, namely $x_1,\ldots,x_8\in\mathbb{Z}$ and $x_9 \in \mathbb{N}$,
and thus finally obtained $\nu \leq 11$ where the unknowns range over all integers.
For this reason, we focus our work on Sun's result, even if the condition $\nu \leq 9$
seems better than $\nu \leq 11$.
This is also our main justification for the work to construct $J_{1+q,\mathbb{C}}$
with an arbitrary $q$.
 We use it for
$q=3$, $q=7$ as required in \cite{MR75},
and $q=17$ in Sun's number theoretic results \cite{sun2021results}.

\section{Preliminaries}\label{s:det}

We shortly remind the definition of {\it diophantine sets}, that will be used
in the formalization. A diophantine polynomial in $k$ variables $v_1, v_2, \ldots, v_k$ is a
linear combination of monomials with non-zero coefficients of the shape $c\cdot v_1^{p_1}v_2^{p_2}v_3^{p_3} \ldots v_j^{p_j}$, where the coefficients $c$ are
integers, the exponents
$p_i$ are natural numbers, and $v_i$ are variables.
The variables will be separated into parameters and unknowns as follows.
A diophantine equation, is an equation of the form $P(x_1, \ldots, x_j, y_1,\ldots, y_k) = 0$, where $P$ is a diophantine polynomial and
$x_1,\ldots, x_j$, $y_1,\ldots, y_k$ indicate the parameters and unknowns, respectively.
A set $D\subseteq \mathbb{N}^n$ of $n$-tuples is called diophantine if there exists a $n + k$-variable diophantine polynomial $P$
such that \mbox{$\langle x_1,\ldots,x_n\rangle\!\in\!D$} if and only if
there exist unknowns $y_1,\ldots,y_k \in \mathbb{N}$ such that $P(x_1, \ldots, x_j, y_1,\ldots, y_k) = 0$.
Similarly, an $n$-ary predicate $\mathcal{P}$ is diophantine, iff the set of $n$-tuples for which the predicate
$\mathcal{P}$ is satisfied is diophantine.
In particular, the divisibility relation $\mathcal{P}(a,b)\!\equiv\! a \mid b$ and
congruence $\mathcal{P}(a,b,c)\!\equiv\!a\mid b$ are diophantine, as
$\mathcal{P}(a,b)\!\iff\!\exists_x\:ax\!-\!b=0$ and $\mathcal{P}(a,b,c)\!\iff\!\exists_x\:a\!-\!b\!-\!cx =0$.

Note, that the number of variables used in a polynomial defining a diophantine property includes
the explicitly stated parameters $a,b,c$ but also the implicitly appearing unknown $x$.
Informally, a function is referred to as  diophantine if the relation between its arguments and its results is.
In particular, the key lemma of Matiyasevich, stating that the exponential function is diophantine means that the
relation $\mathcal{P}(a,b,c)\equiv a = b^c$ is.

\section{Pell Equation}
Even if Matiyasevich originally showed the key lemma using properties
of the Fibonacci sequence,
further results and publications in the domain use the Pell equation instead.
The Pell equation states $x^2\!-\!Dy^2=1$ with a non-square parameter $D$.
If $D=a^2\!-\!1$, we can explicitly give all
the solutions of this equation via Lucas sequences:
$x=\chi_a(n), y=\psi_a(n)$:
\begin{align*}
\chi_a(0)&=1,&\chi_a(1)&=a,&\chi_a(n+2)&=2a\chi_a(n+1)-\psi_a(n),\\
\psi_a(0)&=0,&\psi_a(1)&=1,&\psi_a(n+2)&=2a\psi_a(n+1)-\psi_a(n).
\end{align*}
which is the approach used in
the HOL-Light, 
the Lean \cite{carneiro2018matiyasevic} and Mizar \cite{hilb10_1} formalizations.
Alternatively, one can consider the equation $(ax-y)^2 - (a^2-1)x^2 =1$.
This, transformed as $x^2 - bxy+y^2=1$ with $b=2a$ can be used to build the sequence
of solutions $\alpha_b(n)$. That sequence has the interesting property:
if $x^2 - bxy+y^2=1$ then either $x=\alpha_b(m), y=\alpha_b(m+1)$ or
$x=\alpha_b(m+1), y=\alpha_b(m)$. The latter approach is used in the Coq \cite{H10coq} and
Isabelle \cite{DBLP:conf/itp/BayerDPSS19} formalizations. The similarity can be analysed by noticing the relation
$\alpha_b(n)=\psi_{2a}(n)$. Of course $\psi$ is more general, while $\alpha$ has more properties.
However, the condition
$\alpha_b(k)\mid \alpha_n(m) \Leftrightarrow k\mid m$
  (see Equation (3.23) in \cite{Mat2000})
 is explicitly stated in the publication describing the Coq formalization \cite{H10coq} and the
 previous version of the Isabelle formalization \cite{DBLP:conf/itp/BayerDPSS19}
while $\psi_a(k)\mid \psi_n(m) \Leftrightarrow k\mid m$ is
proved in HOL Light, Lean \cite{carneiro2018matiyasevic} and in Mizar~\cite{hilb10_1}\footnote{
See the complete formalization statement of
the theorem \texttt{Y\_DIVIDES} in HOL-Light
\url{https://github.com/jrh13/hol-light/blob/master/Examples/pell.ml},
the theorem
\texttt{y\_dvd\_iff}
in Lean \url{https://github.com/leanprover-community/mathlib/blob/master/src/number_theory/pell.lean},
theorems \texttt{HILB10\_1:34} and \texttt{HILB10\_1:36}
in Mizar
  \url{http://mizar.uwb.edu.pl/version/current/html/hilb10_1.html}.%
%
%
}.

\newcommand\EBox[0]{=\mkern-6mu\Box}

Irrespective of the considered sequence, all formalizations prove that
$a = \alpha_b(c)$ or $a = \psi_b(c)$
can be represented using an (implicit) diophantine relation, stated as a combination
of less complicated diophantine relations. Additionally these sub-relation
use additional explicit unknowns and may also have implicit ones in these relations.
For example Bayer et al. \cite{DBLP:conf/itp/BayerDPSS19} express $3 < b \wedge a = \alpha_b(c)$ using
6 explicit unknowns and 15 relations including, e.g., 4 uses of equivalence $\equiv$.
Carneiro \cite{carneiro2018matiyasevic} uses 5 additional unknowns explicitly
and several congruences.
Similarly, our previous work used 6 explicitly given additional unknowns \cite{hilb10_1},
reduced to 5 explicit unknowns and a single implicit ones \cite{hilb10_6} in order to
achieve the 26-variable polynomial \eqref{poly26}.
Here, in order to formalize the best known 10-variable polynomial proposed in \cite{eng_30}, we use
the representation proposed by Matiyasevich and Robinson \cite{MR75}:
Two explicit unknowns $i,j$ and 3 implicit ones corresponding to the relations
$\EBox$ ($\EBox$ is a one-argument relation, which is true when the argument is a square of a natural
number, i.e., $x=\Box\Leftrightarrow\exists_{n\in \mathbb{N}} x=n^2$), $\mid$ and $\leq$.

Note that the \autoref{PellSystem} also depends on the parameter $e \in\mathbb{N}$,
which can be simply eliminated by $e=0$. However, we use the original formulation
to simplify the comparison of conditions in \autoref{PellSystem} and \autoref{total}, where we
substitute $e=L\!-\!1$.


\begin{theorem}\label{PellSystem} (\verb|HILB10_8:19|)
Let $A,B,C\in\mathbb{N}$ with $A>1$, $B>0$ and $e \in\mathbb{N}$. Then $C = \psi_A(B)$ if and only if there exists
$i,j\in\mathbb{N}$ and auxiliary unknowns $D,E,F,G,H,I\in \mathbb{Z}$ such that
\begin{equation}\label{PellSystemEQ}
DFI\EBox\: \wedge\: F \mid (H - C) \:\wedge\: B \leq C
\end{equation} and
    $D= (A^2\!-\!1)C^2\!+\!1$,
   $E= 2(i\!+\!1)D(e\!+\!1)C^2$,
   $F= (A^2\!-\!1)E^2+1$,
    $G = A\!+\!F(F\!-\!A)$,
   $H = B\!+\!2jC$,
    $I = (G^2\!-\!1)H^2\!+\!1$, where the auxiliary unknowns can be replaced by polynomials over $A,B,C,i,j$ and $e$.
\end{theorem}
Therefore, $C = \psi_A(B)$ can be represented as
$0=(DFI\!-\!\alpha^2)^2 + (F\beta\!-\!H\!+\!C)^2(F\beta\!+\!H\!-\!C)^2+(B\!+\!\gamma\!-\!C)^2$, where
$\alpha,\beta,\gamma\in\mathbb{N}$ are hidden unknowns.

\section{Prime Numbers}

The main idea behind the construction of a polynomial representing the prime numbers uses Wilson's theorem,
i.e., for any positive integer
$k$, $k + 1$ is prime if and only if
$k + 1\mid k! + 1$.
Note that in Mizar we had to formalize the fact that $y = x!$ is a diophantine relation,
since it is one of the key steps to proving the DPRM-theorem in our approach.
A proof that only focuses on the existence of this polynomial
can be expressed in a surprisingly concise way (less than 100 Mizar lines of proof) using
higher-order schemes. Compare this with more than 2000 lines required to prove
that for any $k\in\mathbb{N}^+$ holds $k + 1$ is prime if and only if,
there exist $a-z\in\mathbb{N}$ unknowns for which
\begin{multline}\label{poly26}
(wz\!+\!h\!+\!j\!-\!q)^2 + ((gk\!+\!g\!+\!k)(h\!+\!j)\!+\!h\!-\!z)^2 +
        ((2k)^3(2k\!+\!2)(n\!+\!1)^2+1\!-\!f^2)^2+\\
        (p\!+\!q\!+\!z\!+\!2n\!-\!e)^2 +
        (e^3(e\!+\!2)(a\!+\!1)^2\!+\!1\!-\!o^2)^2 +
        (x^2\!-\!(a^2\!-\!1)y^2\!-\!1)^2 +
        (16(a^2\!-\!1)r^2y^2y^2\!+\!1 -u^2)^2+\\
        (((a\!+\!u^2(u^2\!-\!a))^2\!-\!1)(n\!+\!4dy)^2\!+\!1\!-\!(x\!+\!cu)^2)^2+
        (m^2\!-\!(a^2\!-\!1)l^2\!-\!1)^2 + (k\!+\!i(a\!-\!1)\!-\!l)^2 +\\
        (n\!+\!l\!+\!v\!-\!y)^2 +(p\!+\!l(a\!-\!n\!-\!1)+b(2a(n\!+\!1)\!-\!(n\!+\!1)^2\!-\!1)\!-\!m)^2+\\
        (q\!+\!y(a\!-\!p\!-\!1)\!+\!s(2a(p\!+\!1)\!-\!(p\!+\!1)^2\!-\!1)\!-\!x)^2+
        (z\!+\!pl(a\!-\!p)\!+\!t(2ap\!-\!p^2\!-\!1)-pm)^2
\end{multline}
equals zero.

To get closer to the 10 variables, we define the notion of prime numbers following \cite{eng_30}:

\begin{theorem}\label{total}(\verb|HILB10_8:23|)
Let $k \in \mathbb{N}$. Then $k$ is prime if and only if there exists
$f,i,j,m,u\in \mathbb{N}^+$, $r,s,t\in\mathbb{N}$ unknowns and auxiliary unknowns
        $A-I,L,M,S-W,Q\in \mathbb{Z}$ such that
\begin{eqnarray}\nonumber
&&DFI \EBox \wedge
      (M^2\!-\!1)S^2\!+\!1 \EBox \wedge
      ((MU)^2 -1)T^2\!+\!1 \EBox \wedge\\
\nonumber
&&   (4f^2 -1)(r-mSTU)^2 + 4u^2S^2T^2 < 8fuST(r-mSTU) \wedge
\\
 &&FL \mid (H-C)Z + F(f+1)Q + F(k+1)((W^2-1)Su-W^2u^2 +1)\label{PrimeEQ3}
\end{eqnarray}
and $A = M(U\!+\!1)$,
$B = W\!+\!1$,
$C = r\!+\!W\!+\!1$,
$D= (A^2\!-\!1)C^2\!+\!1 $,
 $     E= 2iC^2LD$,
  $    F= (A^2\!-\!1)E^2\!+\!1 $,
   $   G = A\!+\!F(F\!-\!A)$,
    $  H = B\!+\!2(j\!-\!1)C$,
     $ I = (G^2\!-\!1)H^2\!+\!1 $,
      $W = 100fk(k\!+\!1) $,
      $U = 100u^3W^3\!+\!1$,
      $M = 100mUW\!+\!1 $,
      $S = (M\!-\!1)s\!+\!k\!+\!1 $,
      $T = (MU\!-\!1)t\!+\!W\!-\!k\!+\!1 $,
      $Q = 2MW\!-\!W^2\!-\!1 $,
        $L = (k\!+\!1)Q $.
\end{theorem}

One can verify, that the simplest polynomial specified by \autoref{PrimeEQ3},
uses 8 unknowns explicitly along
with one implicit one for each relation (three occurrences of $\EBox$, inequality,
 and divisibility). Together with
$k$ this gives a total of 14 variables.
The fact, that this does not require 5 implicit unknowns, was a striking solution
proposed in \cite{eng_30}.
 In the next section, we will reduce the 5
implicit unknowns used for \autoref{PrimeEQ3} to a single one.

\section{The Polynomial Reduced to a Single Unknown Variable}\label{Poly}

We first notice $a\EBox \equiv \exists_{x\in\mathbb{N}}\: 0 = x^2\!-\!a =  \prod (x\pm\sqrt{a})$,
where the product considers all sign combinations.
Of course, the product is a polynomial, but its factors are not. Additionally, $a$ can be negative,
so we need to consider $\mathbb{C}$ as the domain and take into account the non-uniqueness of the square root
(however, since $a\in\mathbb{Z}$, there is only one root in the first quadrant of the complex plane).
\begin{theorem}\label{MRthm}
Suppose $A_1,\ldots,A_q\in\mathbb{Z}$. Then $A_1\EBox, \ldots, A_q\EBox$
if and only if
$$0 = \prod(X\pm\sqrt{A_1}\pm \sqrt{A_2}W\pm \ldots \pm\sqrt{A_q}W^{q-1})$$
for some $X\in \mathbb{Z}$ where $W=1+A_1^2+\ldots+A_q^2$.
\end{theorem}

\autoref{MRthm} is formulated in \cite{MR75} and used for $q=7$,
however, the formalization additionally requires a justification that the product over $2^q$ possible combinations of signs
eliminates similar elements giving
a linear combination of $\binom{2^{q-1}\!+\!q\!-\!1}{q\!-\!1}$ monomials with non-zero coefficients.
For this, we will define a helper polynomial $J_{n,\mathcal{R}}$ in \autoref{Jdef}. There,
all factors that included square roots will appear in even powers, which will eliminate these roots.
This allows using $J_{q+1,\mathbb{C}}(r_0,\sqrt{r_1W^2},\ldots,\sqrt{r_{n-1}W^{2q\!-\!2}})$
as an appropriate $\mathbb{Z}$-valued polynomial over $q+1$, used in \autoref{MRthm}
following Matiyasevich's elegant adaptation in order
to ensure the satisfiability of all 5 predicates from $\autoref{PrimeEQ3}$.
Our current formalization does not include \autoref{MRthm} in its full generality (ongoing
work with most of the needed lemmas complete), as we only use it with $q=3$ and
substitution the constant $W=2$ instead of the polynomial $W= 1+r_1^2+\ldots+r_q^2$.

\begin{theorem}\label{Jdef}(\verb|POLYNOM9:def 10|)
Let $\mathcal{R}$ be a commutative ring, $n\in \mathbb{N}$ with $n>1$.
There exists an $\mathcal{R}$-valued polynomial over $n$ variables $J_{n,\mathcal{R}}$ obeying the following conditions:
\begin{itemize}
\item $J_{n,\mathcal{R}}(r_1,r_2,\ldots,r_n) = \prod(r_1\pm r_2 \pm \ldots\pm r_{n})$ for all $r_1,r_2,\ldots,r_{n}\in\mathcal{R}$,
\item
let $p_{\alpha}R^\alpha$ be any monomial with nonzero coefficients of $J_{n,\mathcal{R}}$, where
$R^\alpha = r_1^{\alpha_1}\cdot r_2^{\alpha_2} \cdot \ldots \cdot r_n^{\alpha_n}$. Then
every power of $\alpha_i$ is even,
the sum of the factors $\sum_{i=1}^n \alpha_i$ is equal to  $2^{n-1}$, 
coefficient $p_{\alpha}$ is an integer multiple of $1_\mathcal{R}$, i.e., is equal to $1_\mathcal{R}\!+\!1_\mathcal{R}\!+\ldots+\!1_\mathcal{R}$
or $-1_\mathcal{R}\!-\!1_\mathcal{R}\!\ldots-\!1_\mathcal{R}$
and the coefficient of $r_1^{2^{n-1}}$ equals $1_\mathcal{R}$.
\end{itemize}
\end{theorem}
The existence of a polynomial that has these properties is quite an involved proof by induction. We only show here the
outline of the most difficult part. By the induction hypothesis:
$\prod(r_1 \pm \ldots\pm r_n \pm r_{n+1})=$%
$\prod(r_1 \pm \ldots \pm (r_n +r_{n+1}))\;\cdot\;\prod(r_1\pm \ldots \pm (r_n\!-\!r_{n+1}))=$
$\sum_\alpha c_\alpha R^\alpha (r_n\!+\!r_{n+1})^{2i_\alpha}\:\cdot\:\sum_\beta c_\beta R^\beta (r_n\!-\!r_{n+1})^{2i_\beta}$
where $R^\alpha$ represent products of $r_1,\ldots,r_{n-1}$ to even powers.
We multiply the sums as follows.
If $i_\alpha = i_\beta$, then $c_\alpha R^\alpha (r_n\!+\!r_{n+1})^{2i_\alpha}c_\beta R^\beta (r_n\!-\!r_{n+1})^{2i_\beta}=
  c_\alpha c_\beta R^\alpha R^\beta (r_n^2\!-\!r_{n+1}^2)^{i_\alpha}$. If $i_\alpha < i_\beta$ (the case $i_\alpha > i_\beta$ is similar)
we add each two summands
\begin{multline}  c_\alpha R^\alpha (r_n\!+\!r_{n+1})^{2i_\alpha}c_\beta R^\beta (r_n\!-\!r_{n+1})^{2i_\beta}+
   c_\beta R^\beta (r_n\!+\!r_{n+1})^{2i_\beta}c_\alpha R^\alpha (r_n\!+\!r_{n+1})^{2i_\alpha} = \\
  c_\alpha c_\beta R^\alpha R^\beta \cdot (r_n^2\!-\!r_{n+1}^2)^{i_\alpha} \cdot\sum_{i=0}^{i_\beta\!-\!i_\alpha}\,2\cdot
\binom{2(i_\beta\!-\!i_\alpha)}{2i} \,r_n^{2i}\,r_{n+1}^{2(i_\beta\!-\!i_\alpha\!-\!i)}.
\end{multline}
In both cases, we obtain a polynomial, where all variables are raised to even powers, which completes the most involved
part of the proof. \hfill $\blacksquare$\\
The complete formalized proof of this theorem includes all the required sign combinations
and required 143 helper lemmas and 13K lines of proofs.


The next step in the simplification of the polynomial given in \autoref{MRthm} (following \cite{eng_30}), proceeds by
defining $K_1(y,x_1,x_2,x_3)$ to be $J_{4,\mathbb{C}}(-y,\sqrt{x_1},\sqrt{4x_2},\sqrt{16x_3})$ and proving
$\exists_{y\in\mathbb{Z}}K_1(y,x_1,x_2,x_3)=0 \Leftrightarrow x_1\EBox\wedge x_2\EBox \wedge x_3\EBox$
under the assumptions $x_1,x_2,x_3\in\mathbb{N}$ and $2\nmid x_1, 2 \nmid x_2$.
Note, that the substitution (compare \autoref{PrimeEQ3})
\[
  x_1=(M^2\!-\!1)S^2\!+\!1$, $x_2=((MU)^2 -1)T^2\!+\!1$, $x_3=DFI
\]
satisfies these assumptions.

The next step in the informal proof performs a rational substitution in the integer polynomial and justifies that
this is again an integer polynomial. This requires some work with the type system in the formalization. Indeed,
we perform the substitution $y:=y-\frac{r}{p}$ in $K_1$, where $p,r$ are the new variables.
In order to construct the polynomial $K_2(y,x_1,x_2,x_3,p,r)$ to be $p^8\cdot K_1(y-\frac{r}{p},x_1,x_2,x_3)$
for $y,x_1,x_2,x_3,p,r\in\mathbb{R}$ where $p\neq 0$, the formalization is split into two stages.
First, we define $K_2^\prime(y,x_1,x_2,x_3,z) = K_1(y-z,x_1,x_2,x_3)$, where the power $\alpha$ of
$z$ in each monomial in $K_2^\prime$ is $\leq 8$ and replace $z^\alpha$ by $p^{\alpha}r^{8-\alpha}$ obtaining $K_2$.
This way, we obtain a polynomial, where $K_2(y,x_1,x_2,x_3,p,r)=0$ ensures $p|r$ for every $p,r\in\mathbb{N}, p\neq 0$,
but both factors in \autoref{PrimeEQ3} are non-negative and $FL>0$.

The final theorem confirms that $K$ is a polynomial:
\begin{theorem}\label{final} (\verb|POLYNOM9:77|)
Let $x_1,x_2,x_3,p,r,n\in\mathbb{N}$, $v\in\mathbb{Z}$ where $2\nmid x_1$, $2 \nmid x_2$, $p>0$, $ n > \sqrt{x_1}+2\sqrt{x_2}+4\sqrt{x_3}+r$.
Then $\exists_{y\in\mathbb{N}} K(y,x_1,x_2,x_3,p,r,n,v) = 0$ iff
$x_1\EBox \wedge x_2 \EBox \wedge x_3\EBox \wedge p\mid r \wedge 0 \leq v,$
where we have $K(y,x_1,x_2,x_3,p,r,n,v) = K_2(y\!-\!nv,x_1,x_2,x_3,p,r)$.
\end{theorem}
With $K$, we can represent \autoref{PrimeEQ3} using a single implicit unknown and we can represent primes
using the 10-variable polynomial with the following substitutions in $K$:
\begin{eqnarray*}
v&=&8fuST(r\!-\!mSTU) - ((4f^2 \!-\!1)(r\!-\!mSTU)^2\! +\! 4u^2S^2T^2)-1,\\
n&=&MS\!+\!2MUT\!+\!4A^2CEGH+2(HL\!+\!FfQ\!+\!Fk(W^2Su\!+\!W^2u^2)).
\end{eqnarray*}
With the abbreviations expanded, this gives a diophantine
polynomial \mizV[Poly]{Poly} of degree > 6000 over parameter $k$ and unknowns $f,i,j,m,u,r,s,t,y$,
so we present only
the non-expanded version with the following property.
Let \mizV[k]{k}$\:\in\mathbb{N}^+$. Then \mizV[k]{k+1} is prime if and only if there exists
a 10-element vector of natural numbers \mizV[v]{v} such that
the first element is equal to \mizV[k]{k} (\mizV[v,k]{v.1'='k})
and \mizV[v]{v} is a root of \mizV[Poly]{Poly}
(\mizV[Poly,v]{eval(Poly,v)'='0.${\bf F}_\mathbb{R}$}).
In the formal proof, rather than specify the existence of a parameter $k$ and 9 unknowns,
we simplify this by using a 10-element vector with one element equal to \mizV[k]{k}. The final
statement in Mizar is:
\begin{Mizar}{Poly,k,v}
theorem :: POLYNOM9:85
  ex'Poly'be'INT-valued'Polynomial'of'10,${\bf F}_\mathbb{R}$'st
    for'k'be'positive'Nat'holds
      k+1'is'prime'iff'ex'v'being'natural-valued'Function'of'10,${\bf F}_\mathbb{R}$'st
                        v.1'='k'&'eval(Poly,v)'='0.${\bf F}_\mathbb{R}$;
\end{Mizar}%
This is already the minimal polynomial and completes our goal.


\begin{thebibliography}{10}

\bibitem{DBLP:conf/itp/BayerDPSS19}
Jonas Bayer, Marco David, Abhik Pal, Benedikt Stock, and Dierk Schleicher.
\newblock The {DPRM} theorem in {I}sabelle (short paper).
\newblock In John Harrison, John O'Leary, and Andrew Tolmach, editors, {\em
  10th International Conference on Interactive Theorem Proving, {ITP} 2019},
  volume 141 of {\em LIPIcs}, pages 33:1--33:7. Dagstuhl, 2019.
\newblock \href {https://doi.org/10.4230/LIPIcs.ITP.2019.33}
  {\path{doi:10.4230/LIPIcs.ITP.2019.33}}.

\bibitem{carneiro2018matiyasevic}
Mario Carneiro.
\newblock A {L}ean formalization of {M}atiyasevi{\v c}'s theorem, 2018.
\newblock \href {http://arxiv.org/abs/1802.01795} {\path{arXiv:1802.01795}}.

\bibitem{DBLP:journals/jsyml/Davis53}
Martin Davis.
\newblock Arithmetical problems and recursively enumerable predicates.
\newblock {\em J. Symb. Log.}, 18(1):33--41, 1953.
\newblock \href {https://doi.org/10.2307/2266325} {\path{doi:10.2307/2266325}}.

\bibitem{DPR1961}
Martin Davis, Hilary Putnam, and Julia Robinson.
\newblock The decision problem for exponential diophantine equations.
\newblock {\em Annals of Mathematics}, 74:425--436, 1961.
\newblock \href {https://doi.org/10.2307/1970289} {\path{doi:10.2307/1970289}}.

\bibitem{Jones82}
James~P. Jones.
\newblock Universal diophantine equation.
\newblock {\em The Journal of Symbolic Logic}, 45(3):549--571, 1982.

\bibitem{AMM76}
James~P. Jones, Daihachiro Sato, Hideo Wada, and Douglas Wiens.
\newblock Diophantine representation of the set of prime numbers.
\newblock {\em The American Mathematical Monthly}, 83(6):449--464, 1976.

\bibitem{Mat2000}
Michael Lamoureux, editor.
\newblock {\em On Hilbert's Tenth Problem}, volume~1. Pacific Institute for the
  Mathematical Sciences, PIMS Distinguished Chair Lectures, 2000.

\bibitem{H10coq}
Dominique Larchey-Wendling and Yannick Forster.
\newblock {Hilbert's Tenth Problem in Coq}.
\newblock In Herman Geuvers, editor, {\em 4th International Conference on
  Formal Structures for Computation and Deduction (FSCD 2019)}, volume 131 of
  {\em LIPIcs}, pages 27:1--27:20, Dagstuhl, Germany, 2019. Dagstuhl.
\newblock \href {https://doi.org/10.4230/LIPIcs.FSCD.2019.27}
  {\path{doi:10.4230/LIPIcs.FSCD.2019.27}}.

\bibitem{YuriOrg}
Yuri Matiyasevich.
\newblock Enumerable sets are diophantine.
\newblock {\em Doklady Akademii Nauk SSSR (in Russian)}, 191:279--282, 1970.

\bibitem{eng_30}
Yuri Matiyasevich.
\newblock Primes are nonnegative values of a polynomial in 10 variables.
\newblock {\em Journal of Soviet Mathematics}, 15:33--44, 1981.
\newblock \href {https://doi.org/10.1007/BF01404106}
  {\path{doi:10.1007/BF01404106}}.

\bibitem{MR75}
Yuri Matiyasevich and Julia Robinson.
\newblock Reduction of an arbitrary diophantine equation to one in 13 unknowns.
\newblock {\em Acta Arithmetica}, 27:521--553, 1975.

\bibitem{hilb10_5}
Karol P{\k{a}}k.
\newblock {F}ormalization of the {MRDP} theorem in the {M}izar system.
\newblock {\em Formalized Mathematics}, 27({\bf 2}):209--221, 2019.
\newblock \href {https://doi.org/10.2478/forma-2019-0020}
  {\path{doi:10.2478/forma-2019-0020}}.

\bibitem{hilb10_6}
Karol P{\k{a}}k.
\newblock Formalization of prime representing polynomial in {M}izar.
\newblock In {\em FMM 2021 workshop}, 2021.
\newblock URL:
  \url{http://alioth.uwb.edu.pl/~pakkarol/articles/KP-FMM2021.pdf}.

\bibitem{hilb10_1}
Karol P\k{a}k.
\newblock The {M}atiyasevich {T}heorem. {P}reliminaries.
\newblock {\em Formalized Mathematics}, 25({\bf 4}):315--325, 2017.
\newblock \href {https://doi.org/10.1515/forma-2017-0029}
  {\path{doi:10.1515/forma-2017-0029}}.

\bibitem{JR69}
Julia Robinson.
\newblock Diophantine decision problems.
\newblock {\em Studies in number theory}, 6:76--116, 1969.

\bibitem{LNT:Smorynski}
Craig~Alan Smorynski.
\newblock {\em Logical Number Theory {I}, An Introduction}.
\newblock Universitext. Springer-Verlag Berlin Heidelberg, 1991.

\bibitem{sun2021results}
Zhi-Wei Sun.
\newblock Further results on {H}ilbert's {T}enth {P}roblem.
\newblock {\em Science China Mathematics}, 64:281--306, 2021.
\newblock \href {https://doi.org/10.1007/s11425-020-1813-5}
  {\path{doi:10.1007/s11425-020-1813-5}}.

\end{thebibliography}

\end{document}